\documentclass[12pt]{amsart}

\usepackage{latexsym}
\usepackage{amsmath}
\usepackage{amssymb}
\usepackage{amsfonts}
\usepackage{amsthm}
\usepackage{amsrefs}
\usepackage{mathrsfs}

\newcommand {\Z}  {\mathbb{Z}}

\def\newNumber#1{%
\refstepcounter{thm}\label{#1}}

\newcommand{\namedRef}[1]{%
\csname env:#1\endcsname\ \ref{#1}%
}

\def\envLRT#1#2{\expandafter\gdef\csname env:#1\endcsname{#2}}

\makeatletter
\newcommand{\namedLabel}[1]{\expandafter\xdef\csname env:#1\endcsname{\envName}%
\immediate\write\@auxout{\string\envLRT{#1}{\envName}}%
\label{#1}}
\newcommand{\NamedLabel}[1]{\expandafter\xdef\csname env:#1\endcsname{\envName}%
\immediate\write\@auxout{\string\envLRT{#1}{\envName}}%
\label{#1}}
\makeatother

\def\newType#1{\newtheorem{#1}[thm]{#1\gdef\envName{#1}}}
\newtheorem{thm}{Theorem}[section]
\newtheorem{Theorem}[thm]{Theorem\gdef\envName{Theorem}}
\newType{Lemma}
\newType{Corollary}
\newType{Proposition}

\theoremstyle{definition}
\newtheorem{Definition}[thm]{Definition\gdef\envName{Definition}}
\newtheorem{Remark}[thm]{Remark\gdef\envName{Remark}}
\newtheorem{Remarks}[thm]{Remarks\gdef\envName{Remarks}}

\newcommand{\cy}[1]{\Z/{#1}\Z}

\newcommand{\rightlabeledarrow}[3][]{%
\def\LRTxx{#1}\ifx\LRTxx\empty\relax%
\setbox0=\hbox{$\mathop{\expandafter\hbox%
{\rightarrowfill}}\limits_{\hbox{\ $\scriptstyle#3$\ }}^{\hbox{\ $\scriptstyle#2$\ }}$}%
\dimen0=\wd0\advance\dimen0 by 4pt%
\edef\LRTxx{\the\dimen0}\else%
\setbox0=\hbox to \LRTxx{$\mathop{\expandafter\hbox%
{\rightarrowfill}}\limits_{\hbox{$\scriptstyle#3$}}^{\hbox{$\scriptstyle#2$}}$}%
\dimen0=\wd0\advance\dimen0 by 4pt%
\fi%
\setbox0=\hbox{$\mathop{\expandafter\hbox%
{\rightarrowfill}}\limits_{\hbox{$\scriptstyle#3$}}^{\hbox{$\scriptstyle#2$}}$}%
\ifdim\dimen0<\wd0 \relax\dimen0=\wd0 \advance\dimen0 by 4pt\fi%
\ \hbox to\dimen0{\hfill\hbox to 0pt{\hss$\mathop{\hbox to\LRTxx {\rightarrowfill}}%
\limits_{\hbox to 0pt{\hss{$\scriptstyle#3$}\hss}}^{\hbox to 0pt{\hss{$\scriptstyle#2$}\hss}}$\hss}\hfill}\ %
}

\newcommand{\downlabeledarrow}[3][]{%
\raise2pt\hbox to 0pt{\hss$\scriptstyle #2$}%
\def\LRTxx{#1}\ifx\LRTxx\empty\relax\downarrow%
\else #1\downarrow\fi\raise2pt\hbox to 0pt{$\scriptstyle #3$\hss}%
}
\newcommand{\downlabeledarrowadjusted}[5][]{%
\raise#4\hbox to 0pt{\hss$\scriptstyle #2$}%
\def\LRTxx{#1}\ifx\LRTxx\empty\relax\downarrow%
\else #1\downarrow\fi\raise#5\hbox to 0pt{$\scriptstyle #3$\hss}%
}

\begin{document}

\def\tts{\ifmmode {\mathcal Z}^t \else ${\mathcal Z}^t$\fi}
\def\zts{\Lambda}
\newcommand{\sw}[1]{\def\xx{#1}\ifx\xx\empty\relax w_1 \else w_1\left(#1\right)\fi}
\newcommand{\bo}{BO}
\newcommand{\bso}{BSO}
\newcommand{\obj}{{\rm Obj}}
\def\ttstop{\tts(\Top)}
\def\subtop{\tts(\mathcal T\,)}
\def\objsubtop{\obj\left(\subtop\right)}
\def\objtts{\obj\left(\ttstop\right)}
\def\cohtwist{\omega}
\def\hzero{\mathfrak c}
\def\ztw{$\Z^t$-}
\def\rp#1{\mathbb R\mathbb P^{#1}}
\def\pol#1{s_{#1}}
\def\extendpol#1#2{s_{{\scriptscriptstyle#2}, #1}}
\def\Top{Top}
\def\trivzts{\theta}
\def\segal#1{B#1}
\def\bsegal#1#2{{\bf B}(#1,#2)}
\def\epg#1{\Pi\left(#1\right)}
\def\basept#1{b_{#1}}
\def\cc#1{c_{#1}}

\title{Unoriented geometric functors}
\author{Laurence R. Taylor}
\email{taylor.2@nd.edu}
\thanks{Partially supported by the N.S.F.}

\begin{abstract}
Farrell and Hsiang \cite{FH}*{p.~102} noticed that \cite{LRT} implies that the geometric surgery groups
defined in \cite{Wall}*{Chapter 9} do not have the naturality Wall claims for them. 
Augmenting Wall's definitions using spaces over $\rp{\infty}$ and 
line bundles they fixed the problem.

The definition of geometric Wall groups involves homology with local coefficients and 
these also lack Wall's claimed naturality. 

One would hope that a geometric bordism theory involving non-orientable 
manifolds would enjoy the same naturality as that enjoyed by homology with
local $\Z$ coefficients. 
A setting for this naturality entirely in terms of local $\Z$ coefficients is presented in this paper.

Applying this theory to the example of non-orientable Wall groups restores 
much of the elegance of Wall's original approach. 
Furthermore, a geometric determination of the map induced by conjugation by 
a group element is given as well as a discussion of further cases beyond the reach of \cite{LRT}.
\end{abstract}
\maketitle

\section{A review of local coefficients}
A \emph{local coefficient system} is a
functor from the path groupoid, $\epg{X}$, of a space $X$ to some
category, \cite{Spanier}*{page 58}.
Two coefficient systems are \emph{equivalent}
if there is a natural transformation between the two functors.
In the case of interest here the category is the group $\Z$ and 
such a local coefficient system will be called a \ztw system.
\begin{Definition}
For any space $X$, let  $\tts(X)$ denote the category whose objects are 
\ztw systems on $X$ and whose morphisms are 
the natural transformations between them.
\end{Definition}

The data for such a system on a space $X$ can be packaged as a function
$\zts\colon \epg{X} \to \{\pm1\}$ which is a homomorphism of groupoids. 
One such system is the \emph{trivial} \ztw system which assigns $+1$ to
every path.
A natural transformation, or morphism, or just map, between $\zts_0$ and $\zts_1$ 
is a function $\zeta\colon X \to \{\pm1\}$ such that for every path $\lambda\in \epg{X}$,
$$\zts_0(\lambda)\cdot\zts_1(\lambda) = \zeta\bigl(\lambda(0)\bigr)\cdot \zeta\bigl(\lambda(1)\bigr)~.$$

Given any \ztw system $\zts$ on $X$ 
check that for any based loop $L$ in $X$, $\zts(L)$ is independent of the base point
and defines a homomorphism $H_1(X;\cy{2}) \to \{\pm1\}$.
Define the \emph{twist} of $\zts$ to be the corresponding element 
$\cohtwist(\zts)\in H^1(X;\cy{2})$.
Furthermore, check that two \ztw systems $\zts_0$ and $\zts_1$ on $X$ are
equivalent if and only if $\cohtwist(\zts_0) = \cohtwist(\zts_1)$.
If $\cohtwist(\zts)$ is trivial, $\zts$ is said to be \emph{untwisted}.

\medskip
Here is a description of all \ztw systems on a space $X$.
Pick a set of base points $\{ \basept{i}\}$, one in each path component of $X$.
For each point $x\in X$ there is a unique $i$ such that $x$ and $\basept{i}$ can be joined 
by a path.
For each $x$ pick one path $\lambda_x$ joining $x$ to $\basept{i}$.
Use the constant path to join $\basept{i}$ to $\basept{i}$.
This choice is called a \emph{polarization} of $\epg{X}$ and can be characterized as a 
function $\pol{}\colon X \to \epg{X}$ such that $\pol{}(x)(0) = x$, 
$\pol{}(x)(1) = \basept{i}$ and $\pol{}(\basept{i})$ is constant.
\begin{Lemma}\namedLabel{main lemma}
Fix any $\cohtwist\in H^1(X;\cy{2})$ and any polarization $\pol{}$ of $\epg{X}$.
Then any function
$f\colon X \to \{\pm 1\}$ such that $f(\basept{i}) = 1$ for all $i$ 
can be uniquely extended to a  \ztw system $\zts$ such that $\cohtwist(\zts) = \cohtwist$ 
and $f = \zts \circ \pol{}$. 
\end{Lemma}
\begin{proof}
Given a path $\lambda$ from $x_0$ to $x_1$, both points lie in a path component with base point
$\basept{i}$.
There exists a unique $\gamma \in \pi_1(X, \basept{i})$ such that $\gamma$ is equivalent to
$s(x_0) \cdot \gamma \cdot s(x_1)^{-1}$.
Define $\zts(\lambda) = f(x_0)\cdot \cohtwist(\gamma)\cdot f(x_1)$.
\end{proof}

Given a subspace $A\subset X$ and a polarization $\pol{A}\colon A \to \epg{A}$, 
construct a polarization $\extendpol{s_A}{X}\colon X \to \epg{X}$ as follows.
Pick base points $\{ \basept{j} \} \in X$. 
If $a_i\in A$ is one of the base points in $A$, 
let $\extendpol{s_A}{X}(a_i)$ be any path from $a_i$ to the appropriate $\basept{j}$.
If $x\in A$, define $\extendpol{s_A}{X}(x) = \pol{A}(x) \cdot \extendpol{s_A}{X}(a_i)$ 
where $x$ is joined in $A$ to $a_i$.
If $x\in X - A$ then choose any path from $x$ to the appropriate $\basept{j}$ using the constant path 
whenever possible.

\begin{Proposition}
Let $A\subset X$ and suppose $\zts_A$ is a \ztw system on $A$. 
If there exists a class $\cohtwist \in H^1(X;\cy{2})$ which restricts to 
$\cohtwist(\zts_A)$, then there exists $\zts_X$ on $X$ which restricts to $\zts_A$ 
and which satisfies $\cohtwist(\zts_X) = \cohtwist$.
\end{Proposition}
\begin{proof}
Fix a polarization $\pol{A}$ of $A$ and a polarization $\extendpol{s_A}{X}$ of $X$.
Define $f\colon X \to \{\pm 1\}$ as follows.
If $x\in A$, define $f(x) = \zts_A\bigl( \pol{A}(x)\bigr)$ and if $x\in X - A$ let $f(x) = 1$.
Let $\zts_X$ be the \ztw system constructed by \namedRef{main lemma} using $\cohtwist$ 
for the twist.
\end{proof}

Note that if $\zeta$ is a morphism, so is $-\zeta$.
Indeed the sign can be switched or not on each path component so 
the set of morphisms between two \ztw systems on $X$ is an $H^0(X;\cy{2})$-torsor.
An element $\hzero\in H^0(X;\cy{2})$ is equivalent to a homomorphism $H_0(X;\cy{2}) \to \{\pm 1\}$
and hence to a function $X\to \{\pm 1\}$ which is constant on path components.
If $\zeta$ is a morphism, let $\hzero\bullet \zeta$ be defined by 
$\hzero\bullet \zeta(x) = \hzero(x)\cdot \zeta(x)$ for all $x\in X$.

\begin{Proposition}
Let $\iota \colon A\subset X$ and suppose $\zts_0$ and $\zts_1$ are \ztw systems on $X$ with
$\cohtwist(\zts_0) = \cohtwist(\zts_1)$. 
Suppose $\zeta_A\colon \iota^\ast(\zts_0) \to \iota^\ast(\zts_1)$ is given.
Then $\zeta_A$ extends to $\zeta_X \colon \zts_0 \to \zts_1$ provided one of the 
properties below holds. 
\begin{enumerate}
\item $H^0(X;\cy{2}) \rightlabeledarrow{\iota^\ast}{} H^0(A;\cy{2})$ is onto. 
\item Pick base points $\{ a_i \}$ for $A$ and for each path component of $X$ pick 
a base point $\basept{j}$ to be one of the $a_i$.
Pick paths $\lambda_i$ from $a_i$ to the appropriate $\basept{j}$.
For each $i$ require $\zts_0(\lambda_{i}) \cdot \zts_1(\lambda_{i}) = 
\zeta_A\bigl(\lambda_{i}(0)\bigr) \cdot \zeta_A\bigl(\lambda_{i}(1)\bigr)$.
\end{enumerate}
\end{Proposition}
\begin{proof}
Since the twists are the same, there exists $\zeta^\prime_X\colon \zts_0 \to \zts_1$ which 
restricts to $\zeta^\prime_A\colon \iota^\ast(\zts_0) \to \iota^\ast(\zts_1)$.
Let $\hzero_A\in H^0(A;\cy{2})$ denote the element such that 
$\zeta_A = \hzero_A \bullet \zeta^\prime_A$.
If (1) holds, pick $\hzero_X\in H^0(X;\cy{2})$ which restricts to $\hzero_A$ and check
$\zeta_X = \hzero_X \bullet \zeta^\prime_X$.

For case (2), pick base points $\{ \basept{j} \}$ for $X$ so that each $\basept{j} = a_{i_j}$.
Since $\cohtwist(\zts_0) = \cohtwist(\zts_1)$ choose an equivalence $\zeta_X\colon \zts_0 \to \zts_1$.
With a bit of care, further require that $\zeta_X( a_{i_j}) = \zeta_A(a_{i_j})$ for all $j$.
It follows from the hypotheses that $\zeta_X\vert_A = \zeta_A$.
\end{proof}

\begin{Proposition}\namedLabel{extend system and isomorphism}
Let $\zts_0$ be a \ztw system on $X$ and let $\zts_A$ be a \ztw system on $A$ where 
$\iota\colon A \to X$ is the inclusion.
Let $\zeta \colon \zts_A \to \iota^\ast(\zts_0)$ be a map. 
Then there exists a \ztw system $\zts_1$ on $X$ such that $\zts_1$ restricted to $A$ is equal 
to $\zts_A$ and $\zeta$ extends to a map $\zeta \colon \zts_1 \to \zts_0$.
\end{Proposition}
\begin{proof}
Pick a polarization $\pol{X}\colon X \to \epg{X}$ and let $f_0\colon X \to \{\pm 1\}$ be 
$\zts_0 \circ \pol{X}$.
Define $f_1\colon X \to \{\pm 1\}$ by $f_1(x) = f_0(x) \cdot \zeta(x)$ if $x\in A$ and $1$ otherwise.
Then $\zts_1$ is the \ztw system from \namedRef{main lemma}.
\end{proof}

\section{The category $\ttstop$}
A map $f\colon X_0\to X_1$ induces a functor $\epg{X_0} \rightlabeledarrow{}{} \epg{X_1}$ 
and hence a functor $f^\ast\colon \tts(X_1) \to \tts(X_0)$.
Composition and identity behave correctly, so there is
a functor from the category of topological spaces and
continuous functions to the category of categories
which takes $X$ to the category $\tts(X)$ and $f$
to the functor $f^\ast$.

Define a category $\ttstop$ whose objects are all pairs $(X, \zts_X)$ 
for $X$ a topological space and $\zts_X \in \tts(X)$.
A morphism $(X_0, \zts_{X_0}) \to (X_1, \zts_{X_1})$ is a pair $( f, \zeta)$
where $f\colon X_0 \to X_1$ is any map 
and $\zeta\colon \zts_{X_0} \to f^\ast(\zts_{X_1})$ is any natural transformation.

Composition is given by $(g, \psi) \circ (f,\zeta) = ( g \circ f, f^\ast(\psi)\circ \zeta)$.
The identity for $(X,\zts)$ is given by $(1_X,1_{\zts})$.
\bigskip

\begin{Definition}
The category $\ttstop$ has a pairing resembling a product.
Given $(X_i, \zts_{X_i}) \in \objtts$, $i=0$, $1$, define a new object
$$(X_0\times X_1, \zts_{X_0}\bullet\zts_{X_1})~.$$
Recall $\epg{X_0\times X_1} = \epg{X_0}\times \epg{X_1}$ so define 
$$\zts_{X_0}\bullet\zts_{X_1}(\lambda_0\times \lambda_1) = 
\zts_{X_0}(\lambda_0)\cdot \zts_{X_1}(\lambda_1)~.$$
This construction is not a categorical product since there is no projection onto $X_0$
unless $\zts_{X_1}$ is untwisted.
But if $(g_i,\zeta_i)\colon (Y,\zts_Y) \to (X_i,\zts_{X_i})$ are given, there is a morphism
$$(g_0\times g_1, \zeta_0\bullet\zeta_1) \colon
(Y,\zts_Y) \to (X_0\times X_1, \zts_{X_0}\bullet\zts_{X_1})$$
\end{Definition}

\begin{Definition}\namedLabel{preserve involution}
The category $\ttstop$ has an involution: it is the identity on objects and sends 
$(f,\zeta)$ to $(f, -\zeta)$.
\end{Definition}

\begin{Remarks}\namedLabel{bundle remarks}
Given a bundle over $X$ orient each fibre and associate a \ztw system as follows. 
Given a path in $X$, use the homotopy lifting property to get a map from the fibre over
the initial point to the fibre over the terminal point of the path. 
Assign $+1$ to this path if the map is orientation preserving, $-1$ otherwise.
This is a \ztw system and the twist is the first Stiefel-Whitney class of the bundle. 
There is no preferred choice of orientations if the bundle is non-orientable and even
if it is, there are two choices over each path component of $X$. 
Similar remarks work for spherical fibrations. 
\end{Remarks}

\section{Homotopy functors}
The remaining sections describe functors out of subcategories of $\ttstop$ or its opposite.
Many of these functors have a homotopy invariance which is described next.

\begin{Definition}
Let $\subtop$ be a subcategory of $\ttstop$ and 
let $\trivzts_{[0,1]}$ denote the trivial system on $[0,1]$.
Say $\subtop$ \emph{has homotopies} provided the following properties hold.
\begin{enumerate}
\item If $(X,\zts) \in \objsubtop$ then
$$\bigl(X\times [0,1], \zts \bullet \trivzts_{[0,1]}\bigr) \in \objsubtop~.$$
\item For $i=0$ and $i=1$ the morphisms
$$(\iota_i, 1_\zts) \colon (X,\zts) \to \bigl(X\times [0,1], \zts \bullet \trivzts_{[0,1]}\bigr)$$
are in $\subtop$ where $\iota_i$ is the evident inclusion. 
Since $\iota_i^\ast\bigl(\zts \bullet \trivzts_{[0,1]}\bigr)$ is identical to $\zts$, 
$1_\zts$ is a morphism between the relevant \ztw systems.
\item \label{homotopies exist}Let $(f_i, \zeta_i) \colon (X,\zts_X) \to (Y, \zts_Y)$, $i=0$, $1$ 
be two maps in $\subtop$ whose $f_i$ are homotopic. 
For any homotopy $F\colon X\times[0,1] \to Y$ require the existence of a map
$\zeta_F\colon \zts_X \bullet  \trivzts_{[0,1]} \to F^\ast(\zts_Y)$ such that 
$(F , \zeta_F)$ is a map in $\subtop$ and the composition 
$(F,\zeta_F) \circ (\iota_0, 1_{\zeta_X})$ is $(f_0,\zeta_0)$.
\item \label{H0 maps exist}
If $(f , \zeta_i) \colon (X, \zts_X) \to (Y, \zts_Y)$ $i=0$, $1$ are two morphisms with 
the same $f$, then there exists $\hzero \in H^0(X;\cy{2})$ such that $\zeta_0 = \hzero \bullet \zeta_1$.
Require that $(1_X, \hzero \bullet 1_{\zts_X})$ is a map in $\subtop$.
\end{enumerate}
\end{Definition}

\begin{Remarks}
The category $\ttstop$ has homotopies.
The map $\zeta_F$ required by (\ref{homotopies exist}) is unique.
If $r\colon X \times [0,1] \to X$ denotes the evident projection, $r$ is a homotopy from
$1_X$ to itself.
The \ztw system $r^\ast(\zts) = \zts \bullet  \trivzts_{[0,1]}$ so if $\zeta_r$ is 
the identity, $(r,1_{\zts \bullet  \trivzts_{[0,1]}}) \circ (\iota_i, 1_{\zts}) = (1_X, 1_{\zts})$ for 
both $i=0$ and $i=1$.
\end{Remarks}

\begin{Definition}
A functor $\Omega_\ast$ defined on $\subtop$ is a \emph{homotopy functor}
provided $\subtop$ has homotopies and 
$$\Omega_\ast(\iota_0, 1_\zts) \circ \Omega_\ast(r,1_{\zts \bullet  \trivzts_{[0,1]}}) = 
1_{\Omega_\ast(X\times[0,1],\zts \bullet \trivzts_{[0,1]})}$$
There is a similar definition for functors out of $\subtop^{op}$.
\end{Definition}

\begin{Lemma}\namedLabel{homotopy functor}
If $\Omega$ is a homotopy functor
the two maps $$\Omega_\ast(\iota_i,1_\zeta) \colon \Omega_\ast(X,\zts) \to
\Omega_\ast(X\times[0,1] , \zts \bullet \trivzts_{[0,1]})$$ $i=0$, $1$ are equal.
\end{Lemma} 
\begin{proof}
The map $\Omega_\ast(r,1_{\zts \bullet  \trivzts_{[0,1]}})$ is an isomorphism inverse to both the maps 
$\Omega_\ast(\iota_0, 1_\zts)$ and $\Omega_\ast(\iota_1, 1_\zts)$.
By uniqueness of inverse, $\Omega_\ast(\iota_0, 1_\zts) = \Omega_\ast(\iota_1, 1_\zts)$. 
\end{proof}

\medskip

Even for a homotopy functor homotopic maps need not induce the same map.
To see the problem, let $F\colon X\times [0,1] \to Y$ be a homotopy between $f_i\colon X \to Y$
so $f_i = F\circ \iota_i$.
Let $\zeta_F\colon \zts_X \bullet \trivzts_{[0,1]} \to F^\ast(\zts_Y)$ be the map such that 
$(F,\zeta_F) \circ (\iota_0, 1_{\zts_X})$ is $(f_0,\zeta_0)$.
Then let $(F,\zeta_F) \circ (\iota_1, 1_{\zts_X}) = (f_1,\zeta^\prime)$ so by \ref{homotopy functor} 
$\Omega_\ast(f_0,\zeta_0) = \Omega_\ast(f_1,\zeta^\prime)$. 

The problem is that $\zeta^\prime$ may not be $\zeta_1$.
Since $\zeta^\prime$ and $\zeta_1$ are maps between the same two \ztw systems, there exists 
a $\hzero_F\in H^0(X ;\cy{2})$ such that $\zeta^\prime = \hzero_F \bullet \zeta_1$. 
Then by (\ref{H0 maps exist}),
$$\Omega_\ast(f_0,\zeta_0) = \Omega_\ast(f_1,\zeta_1) \circ 
\Omega_\ast(1_{X}, \hzero_F\bullet 1_{\zts_X})~.$$
It remains to identify $\hzero_F$.

\begin{Theorem}\namedLabel{compute c_F} For all $x\in X$
$$\hzero_F(x) = \zts_Y\bigl(F(x\times t)\bigr)\cdot \zeta_0(x) \cdot \zeta_1(x)$$
\end{Theorem}
\begin{proof}
By definition $\hzero_F(x) = \zeta_1(x) \cdot \zeta^\prime(x)$.
Moreover $\zeta^\prime(x) = \zeta_F(x\times 1)$ and $\zeta_0(x) = \zeta_F(x\times 0)$.
Also by definition,
$$\zeta_F( x \times 0) \cdot \zeta_F( x \times 1) = 
F^\ast(\zts_Y)(x\times t) \cdot (\zts_X \bullet \trivzts_{[0,1]})(x\times t)~.$$
But $(\zts_X \bullet \trivzts_{[0,1]})(x\times t)  = r^\ast(\zts_X)(x\times t) = 
\zts_X\bigl(r(x\times t)\bigr) = 1$ since $r(x\times t)$ is a constant path.
Similarly $F^\ast(\zts_Y)(x\times t) = \zts_Y\bigl(F(x\times t)\bigr)$.
\end{proof}
\begin{Corollary}\namedLabel{how homotopies act}
Suppose $X$ is path connected and there exists a point $x\in X$ such that $f_0(x) = f_1(x)$.
Then $F(x\times t)$ is a loop $L$ in $Y$ and $\zts_Y\bigl(F(x\times t)\bigr) = \cohtwist(L)$.
\end{Corollary}

If $X$ is path connected, then $\hzero_F = \pm 1$ and 
$\Omega_\ast(1_{X}, \hzero_F\bullet 1_{\zts_X})$ is the identity if $\hzero_F = +1$ and is 
the map induced by the involution on $\ttstop$ if $\hzero_F = -1$ (\ref{preserve involution}).

\begin{Definition}
Any additive category $\mathscr A$ has an involution which is the identity on objects and 
sends $f$ to $-f$ in the abelian group of morphisms between any two objects.
If $\Omega_\ast$ takes values in $\mathscr A$, say $\Omega_\ast$ \emph{preserves the 
involution} provided $\Omega_\ast(f, -\zeta) = - \Omega_\ast(f,\zeta)$, provided both 
$(f,\pm\zeta) \in \subtop$. 
There is a similar definition if $\Omega_\ast$ is a functor out of $\subtop^{op}$.
\end{Definition}

\section{Twisted homology and twisted cohomology}
There is a detailed treatment of this material in \cite{Wh}*{Chapter VI} 
but here is a quick review.
Twisted homology and cohomology come from modifying the
singular chain complex $S_\ast(X)$.
The chain groups are the same, but the boundary is changed.
Recall that $\partial=\sum_i (-1)^i\partial_i$. 
If $\zts$ is a \ztw system,
$\partial^{\zts} =\sum_i (-1)^i\partial_i^{\zts}$ where $\partial_i^{\zts}$ is
defined as follows.
Given a singular $r${--}simplex $\sigma\to X$, $\partial_i\sigma$ is a 
$(r-1)${--}simplex.
Define $\partial_i^{\zts} \sigma$ to be $\partial_i\sigma$ multiplied by
$\zts$ applied to the path obtained by applying $\sigma$ to the
straight line from the barycenter of $\sigma$ to the 
barycenter of $\partial_i\sigma$.
This can be checked to be a chain complex $S^{\zts}_\ast(X)$.
Homology and cohomology with twisted coefficients is defined
as the homology or cohomology of $S^{\zts}_\ast$.

Given a natural transformation $\zeta\colon \zts_0 \to \zts_1$, define
a chain map $S(\zeta)_\ast\colon S^{\zts_0}_\ast(X) \to S^{\zts_1}_\ast(X)$ 
by sending $\sigma$ to $\pm\sigma$ where $\pm1$ is the value of the
natural transformation applied to the barycenter of $\sigma$.
This chain map induces an isomorphism on twisted homology and
cohomology groups.

It follows that $S^\zts_\ast$ and the twisted homology and cohomology  groups are 
functors out of $\ttstop$ or $\ttstop^{op}$.
The twisted homology and cohomology are homotopy functors.
Furthermore $S^\zts_\ast$, $H_\ast$ and $H^\ast$ all preserve involution (\ref{preserve involution}).

\medskip

The Alexander-Whitney diagonal map induces a chain map
$$S^{\zts_{X_0}}_\ast(X_0) \otimes S^{\zts_{X_1}}_\ast(X_1) \to
S^{\zts_{X_0}\bullet\zts_{X_1}}_\ast(X_0\times X_1)$$
so the usual products have twisted versions: given two systems, $\zts_0$ and $\zts_1$,
there is a bilinear cap product
$$\cap \colon H^r(X;\Z^{\zts_0}) \otimes H_{r+s}(X;\Z^{\zts_1}) \to 
H_{s}(X;\Z^{\zts_0\bullet \zts_1})~.$$

\section{Poincar\'e duality spaces}
A remark that the author has found helpful is that 
\emph{non-orientable manifolds can not be oriented}.
Without paying enough attention to naturality one can come away with a vague feeling 
that they can: 
see the remarks at the end of this section.

If $X$ is path-connected, there can be 
at most one class of local coefficients
$\zts$ and one integer $m$ such that $H_m(X;\Z^{\zts}) \cong \Z$ and such that, 
if $[X]$ is a generator,
\newNumber{tpd}
$$\cap [X]\colon H^r(X;\Z^{\zts^\prime}) \to 
H_{m-r}(X;\Z^{\zts\bullet \zts^\prime})\leqno(\ref{tpd})$$
is an isomorphism for all $r$ and all $\zts^\prime$.

If the first Stiefel-Whitney class $\sw{M}\in H^1(M;\cy{2})$ of 
a  connected, closed, compact manifold or Poincar\'e 
duality space is used as the twist, (\ref{tpd}) holds.
When $\sw{M} = 0$ write $H_\ast(M;\Z)$ for homology with trivial \ztw system. 
The choice of a generator  in $H_m(M;\Z)$ 
is often called an orientation and $M$ is said to be oriented.
Given two oriented connected, closed, compact manifolds $M_1^m$ and $M_2^m$ and a map
$f\colon M_1 \to M_2$, $f$  has a degree.

In the non-orientable case the \ztw system is non-trivial and even  
in the orientable case untwisted but  non-trivial systems can appear.
Whenever the \ztw system is non-trivial, the class $[M] \in H_m(M;\Z^{\zts})$ 
will be called a \emph{fundamental class}.
Given a map $f\colon M_0\to M_1$ which preserves the first Stiefel-Whitney class, 
there is an induced map
$$(f, \zeta)_\ast\colon H_m(M_0;\Z^{\zts_{M_0}}) \to H_m(M_1;\Z^{\zts_{M_1}})$$
and $(f, \zeta)_\ast$ has a degree, but if $\zeta$ is replaced by $-\zeta$,
the degree switches sign. 

\begin{Remarks}
The Farrell-Hsiang repair of Wall's problem can be viewed as checking that one can fix 
a \ztw-system over $\rp{\infty}$ and always use the pull-back local system.
Whitehead \cite{Wh}*{Chapter VI} also discusses naturality using pull-back systems. 
This results in the following conundrum.
Consider the projection $p\colon S^{2n} \to \rp{2n}$ and the antipodal map $\alpha$ on $S^{2n}$:
$\alpha$ commutes with $p$ so there is a commutative triangle of twisted homology groups
starting with any $\zts$ on $\rp{2n}$ with non-trivial twist. 
The induced system on $S^{2n}$ is untwisted, but not trivial 
and $\alpha$ ends up having degree $+1$.
This is all correct but unsettling.

In the approach taken here, one can work with the trivial system on $S^{2n}$, in which case 
$\alpha$ has degree $-1$ but now one must choose an isomorphism between the 
system induced by $p$ and the trivial one, say $\zeta$.
Then the projection maps are $(p,\zeta)$ and $(p, -\zeta)$, $\alpha$ has degree $-1$ 
and commutativity still holds.
But this shows that $\rp{2n}$ can not be ``oriented'' by orienting $S^{2n}$ and using
``$p_\ast$'' to pick out a fundamental class.
\end{Remarks}

\begin{Remarks}
Here is another approach to ``orientation''. 
Twisted homology satisfies excision, so
$H_m(M, M - p\, ; \Z^\zts) = H_m( D^m , S^{m-1} ; \Z^\zts \vert_{D^m})$.
Since $D^m$ is simply-connected, $H_m( D^m , S^{m-1} ; \Z^\zts \vert_{D^m}) \cong \Z$. 

There is always a map
$H_m( M ; \Z^{\zts} ) \to H_m(M, M - p\, ; \Z^\zts)$ and if $\cohtwist(\zts) = \sw{M}$, 
a fundamental class $[M] \in H_m( M ; \Z^{\zts} )$ picks out a generator in 
$H_m( D^m , S^{m-1} ; \Z^\zts \vert_{D^m})$.
Pick the unique morphism from the trivial system to $\Z^\zts \vert_{D^m}$ which is $+1$ at 
the point $p$. 
This gives an orientation on $D^m$.

However, whenever $\zts$ is not trivial, given any cover of $M$ by disks there will be disks 
in the cover for which the selected orientation depends on $p$.
\end{Remarks}

\section{Geometric bordism}
As a warm up exercise, consider unoriented bordism.
For the spaces in the category $\subtop$ take $\bo\times [0,1]^k$ for all integers $k\geqslant 0$.
For the \ztw systems take all the ones with non-trivial twist and all their morphisms, 
and for the maps of spaces take all homotopy equivalences.
Fix a \ztw system $\zts_{\bo}$ on $\bo\times [0,1]^k$ with non-trivial twist.
Hence an object in this ``enhanced'' bordism theory is a manifold $M$ with a \ztw 
coefficient system $\zts_M$, a generator $[M]\in H_m(M;\Z^{\zts_M})$, a map
$\nu\colon M \to \bo\times[0,1]^k$ classifying the stable normal bundle 
and a natural isomorphism  $\zeta\colon \zts_M \to \nu^\ast(\zts_{\bo})$.

The definitions for a manifold with boundary present no problems, but  
manifolds should bound which are not actually equal to a boundary but merely ``equivalent'' to one. 
This involves a diffeomorphism of $M$ with an actual boundary and an identification of 
fundamental classes.
Since there is no natural choice of $\zeta$, this requires making 
$(M,\nu, \hzero \bullet [M], \hzero \bullet \zeta)$
equal to $(M,\nu, -[M], -\zeta)$.
Alternatively one can use \namedRef{extend system and isomorphism} to put a \ztw system on $M\times [0,1]$ which alters the sign of $\zeta$ at the two ends.

Bordism defines a homotopy functor on $\subtop$ to abelian groups which preserves orientation.

\begin{Remark}
These definitions can be repeated with $\bso$. 
Oriented bordism has elements of infinite order, so why are not similar elements present in
the unoriented case?
\end{Remark}

To answer this last question, recall that $\bo$ has a homotopy
$F\colon \bo\times [0,1] \to \bo$ which is the identity at both ends and such that for each 
$b \in \bo$, $F(b\times [0,1])$ is a loop which is not null homotopic.
\namedRef{how homotopies act} implies that the identity map has order $2$ 
so every element in the bordism group has order $2$.

In the oriented case no such homotopy exists.

The case of $Pin^{\pm}$ structures is interesting from this point of view as well.
The bundles are not orientable, but $BPin^\pm = BSpin \times K(\cy{2},1)$ so there 
is a homotopy like that for $\bo$.
However the $Pin^{\pm}$ bordism groups are not always of exponent $2$. 
Here the problem is that the homotopy is not a homotopy of lifts over $\bo$.

\section{Geometric surgery}
The definition of Wall's geometric surgery groups 
mimic the unoriented bordism example.
Start with a reference $n$-ad $K$ and fix a \ztw coefficient system, $\zts_K$.
An object in the bordism theory is a Poincar\'e $(n+1)$-ad $X$, a manifold $(n+1)$-ad $M$,
\ztw coefficient systems $\zts_X$ on $X$ and $\zts_M$ on $M$,
classes $[M]\in H_m(M;\Z^{\zts_M})$ and $[X]\in H_m(X;\Z^{\zts_X})$,
maps $g\colon M \to X$ and $f\colon X \to K$, and finally,
isomorphisms $\zeta\colon \zts_X \to f^\ast(\zts_K)$ and $\psi\colon \zts_M \to g^\ast(\zts_X)$
such that $(g,\psi)_\ast$ has degree $1$.
The map $g$ should be covered by a normal map and should have 
some sort of equivalence on part of the boundary depending on 
exactly which variant of the Wall group is being constructed.

It is better to denote the resulting bordism groups by
$L_\ast(K,\zts_K)$ instead of just recording the twist as Wall does.
Since the torsion is irrelevant to these discussions it is suppressed.

With these minor alterations in the definition, the material in Wall's Chapter 9 goes through with
no difficulty, except for one caveat:
the resulting bordism groups are homotopy functors on the full subcategory of $\ttstop$ 
whose spaces are the homotopy type of CW complexes with a finite $2$-skeleton.

Similar results hold for many other variations of surgery theory in the non-orientable case
including the $LS$ groups of \cite{Wall}*{Chapter 11}, the Cappel-Shaneson $\Gamma$ groups
\cite{CS}, and many others.

\section{Surgery groups}
The philosophical import of Wall's Chapter 9 is  that any ``geometric bordism theory'' which
has a $\pi$-$\pi$ theorem depends only on the fundamental group.

More explicitly, if $K$ is path-connected, 
let $u\colon K \to  K(\pi_1(K),1)$ be  a classifying map for the universal cover.
The map $u$ induces an isomorphism 
$$H^1\bigl( K(\pi_1(K) ; \cy{2}\bigr) \to H^1(K;\cy{2})$$
so the set of twists for the two spaces coincide.
There is an explicit model for $K(\pi,1)$ due to Segal \cite{Segal} denoted $\segal{\pi}$ and there 
is a particular \ztw system $\zts_\cohtwist$ with twist $\cohtwist \in H^1(\segal{\pi} ; \cy{2})$.

The explicit model is Segal's classifying space 
of a small category where the category is the group in question.
There is one $0$-simplex which will be the base point.
Each $k$-simplex for $k>0$ is a sequence of $k$ group elements.
If $G$ is the group, let $\segal{G}$ denote the classifying space.
Recall that $\segal{G}$ is a CW complex and it has the homotopy type of a CW complex with 
finite $2$ skeleton if and only if $G$ is finitely generated and finitely presented.

The space $\segal{G}$ has an explicit polarization, $\pol{G}$, defined as follows. 
Each point $x\in \segal{G}$ lies in the interior of a unique simplex and there is always the 
straight line path to the initial vertex of this simplex. 
This gives a path from $x\in \segal{G}$ to the base point. 
Denote the path by $\pol{G}(x)$.

Given a twist $\cohtwist\in H^1(\segal{G} ;\cy{2})$ define $\zts_{\cohtwist}$ to be the \ztw system 
on $\segal{G}$ from \namedRef{main lemma} with polarization $\pol{G}$, 
function $f\colon \segal{G} \to \{\pm1\}$ identically $1$ and twist $\cohtwist$.

If $\psi\colon H \to G$ is a group homomorphism then the induced map 
$\segal{\psi} \colon \segal{H} \to \segal{G}$ is piecewise linear on each simplex
so 
$$\begin{matrix}
\segal{H}&\rightlabeledarrow{\pol{H}}{}& \epg{\segal{H}}\\
\downlabeledarrow{\segal{\psi}}{}&&\downlabeledarrow{}{\epg{\segal{\psi}}}\\
\segal{G}&\rightlabeledarrow{\pol{G}}{}& \epg{\segal{G}}\\
\end{matrix}
$$
commutes. 
It follows that $(\segal{\psi})^\ast ( \zts_{\cohtwist}) = \zts_{\segal{\psi}^\ast(\cohtwist)}$ so
$$(\segal{\psi} , 1_{\zts_{(\segal{\psi})^\ast(\cohtwist)}}) \colon
(\segal{H}, \zts_{(\segal{\psi})^\ast(\cohtwist)}) \to ( \segal{G} , \zts_{\cohtwist}) \in \ttstop~.$$

Since $\zts_{(\segal{\psi})^\ast(\cohtwist)}$ is just the identity map, write
$\bsegal{\psi}{\cohtwist}$ for the map  $(\segal{\psi} , 1_{\zts_{(\segal{\psi})^\ast(\cohtwist)}})$.
Note that for group homomorphisms, $\psi_0\colon G_0 \to G_1$ and
$\psi_1\colon G_1 \to G_2$,
$\bsegal{(\psi_1\circ \psi_0)}{\cohtwist} = \bsegal{\psi_1}{(\segal{\psi_0})^\ast(\cohtwist)}
\circ \bsegal{\psi_0}{\cohtwist}$.

\medskip
Segal's construction has an additional property.
Any $g\in G$ gives an automorphism of $G$, 
$\cc{g}\colon G \to G$ defined by $\cc{g}(h) = g^{-1} h g$.
This inner automorphism is homotopic to the identity since there is a natural transformation 
between the two functors, namely multiplication by $g \in G$.
Let $C_g\colon \segal{G} \times [0,1] \to \segal{G}$ be the homotopy from $1_{\segal{G}}$
at $0$ to $c_g$ at $1$.
The homotopy is not base point preserving unless $g$ is the identity: in fact the path 
$C_g(t)$ is just the $1$-simplex $g$.

\bigskip
\begin{Definition}\namedLabel{depends on pi1}
A \emph{surgery theory} is a homotopy functor $\Omega$ 
on a subcategory $\subtop$ of $\ttstop$ with the following closure property.
Whenever $(K, \zts)$ is an object in $\subtop$, then for any 
$u\colon K \to \segal{\pi_1(K)}$ which induces an isomorphism on $\pi_1$ 
the pair $\bigl(\segal{\pi_1(K)}, \zts_{(u^\ast)^{-1}(\cohtwist(\zts))}\bigr)$ 
is in $\subtop$ 
and for at least one map $\zeta\colon \zts \to u^\ast\bigl( \zts_{(u^\ast)^{-1}(\cohtwist(\zts))}\bigr)$ 
the pair $(u,\zeta)$ is a map in $\subtop$. 
A surgery theory \emph{only depends on $\pi_1$} provided
$$\Omega_\ast(u,\zeta) \colon 
\Omega_\ast (K, \zts ) \to \Omega_\ast\bigl( \segal{\pi_1(K)}, \zts_{(u^\ast)^{-1}(\cohtwist(\zts))}\bigr)$$
an isomorphism.
\end{Definition}

\begin{Remark}
Wall's geometric surgery groups are a surgery theory 
on the full subcategory of $\ttstop$ whose spaces
have the homotopy type of CW complexes with finite $2$-skeleton.  
They preserve the involution. 
If $\ast>4$, then the Wall geometric surgery groups only depend on $\pi_1$.
\end{Remark}

\begin{Theorem}\namedLabel{Taylor's lemma}
Suppose $\Omega_\ast$ is a surgery theory on $\subtop$ and that $G$ is a group 
such that  $(\segal{G},\zts_{\cohtwist})\in \subtop$.
Then $$\Omega_\ast(\bsegal{c_g}{\cohtwist}) \colon \Omega_\ast(\segal{G}, \zts_{\cohtwist}) 
\to \Omega_\ast(\segal{G}, \zts_{\cohtwist}) $$
is the identity if $\cohtwist(g) = +1$ and is 
$\Omega_\ast(1_{\segal{G}}, -1_{\zts_{\cohtwist}})$ if $\cohtwist(g) = -1$.
If $\Omega_\ast$ preserves involution, then $\Omega_\ast(\bsegal{c_g}{\cohtwist})$ 
is minus the identity if $\cohtwist(g) = -1$.
\end{Theorem}
\begin{proof}
Apply \namedRef{compute c_F} and \namedRef{how homotopies act}: note $\zeta_0(x) = \zeta_1(x)$.
\end{proof}
\begin{Remark}
\namedRef{Taylor's lemma} is the analogue of the main result in \cite{LRT} for 
the geometric Wall groups: it implies that the surgery obstruction map \cite{Wall}*{Cor.~9.4.1, p.~90} 
is natural.
\end{Remark}
 \begin{Remark}
This material also implies the geometric version of the conjecture in \cite{LRT} concerning 
the map induced by conjugation.
Without hypotheses the conjecture is false so some explanation may be helpful. 
Consider the pair case, $\iota\colon H\subset G$. 
Then $\segal{H}$ is a subspace of $\segal{G}$ and $\iota^\ast(\zts_{\cohtwist}) = 
\zts_{\iota^\ast(\cohtwist)}$. 
Let 
$\zeta\colon \zts_{\iota^\ast(\cohtwist)} \to \iota^\ast(\zts_{\cohtwist})$ 
be the identity. 
The data for the pair should be written $\Omega_\ast\bigl((\segal{G},\zts_{\cohtwist}),
(\segal{H},\zts_{\iota^\ast(\cohtwist)}), \zeta\bigr)~.$ 

\begin{enumerate}
\item Any homomorphism of pairs $\psi\colon (G_0,H_0) \to (G_1, H_1)$ 
induces a map of long exact sequences.
\item
If $g\in H$ then conjugation by $g$ induces an automorphism of the entire long exact sequence 
of the pair and is multiplication by $\cohtwist(g)$. 
\item
If $g \in G - H$, then $g^{-1} H g$ must be $H$ before $\cc{g}$ is even a map of pairs.
Fix $\psi\colon H \to H$. 
Then there are examples of groups $G$ with elements $g$ such that $g^{-1} H g = H$ and 
the induced map on $H$ is $\psi$.
Hence nothing can be said in general about the map induced by conjugation on  a pair if
$g \notin H$.
\end{enumerate}
\end{Remark}


\end{document}